\documentclass{amsart}
\usepackage{amsmath,amssymb,amsthm}

\usepackage{booktabs}
\usepackage{algorithm}
\usepackage{algpseudocode}
\usepackage{url}

\newtheorem{theorem}{Theorem}[section]
\newtheorem{lemma}[theorem]{Lemma}
\newtheorem{proposition}[theorem]{Proposition}
\newtheorem{corollary}[theorem]{Corollary}
\newtheorem{conjecture}[theorem]{Conjecture}
\theoremstyle{definition}
\newtheorem{definition}[theorem]{Definition}
\newtheorem{heuristic}[theorem]{Heuristic}
\theoremstyle{remark}
\newtheorem{remark}[theorem]{Remark}

\newcommand{\pmx}{p_{\max}}
\newcommand{\lcm}{\operatorname{lcm}}
\newcommand{\Tree}{\mathcal{T}}

\title{The least Carmichael numbers with $k$ prime factors}
\author{Jonathan Webster}
\address{Mathematical Sciences Department, Butler University}
\email{jewebste@butler.edu}
\author{Sylvia Webster}
\address{Mathematical Sciences Department, Butler University}
\date{\today}

\begin{document}

\begin{abstract}
Let $S_k$ denote the least Carmichael number with exactly $k$ prime factors.
Like the later searches of Pinch, our algorithm builds the prime factors of a candidate in increasing order and stops as soon as the preproduct $P$ satisfies $P\lambda(P) > B$.
At that point there is at most one possible cofactor, and we recover all of its prime factors at once in polynomial time.
We prove an unconditional bound on the number of candidates the search examines, in terms of how far $B$ exceeds the product of the first $k$ primes.
When the search is run in doubling passes, this gives a running time of $e^{O(k)} = S_k^{O(1/\log k)}$ under a natural conjecture on the size of $S_k$.
We also show that the constant in the exponent is sharp for the set of square-free numbers that the bound counts.
Finally, we compute $S_k$ for every $k \leq 64$; Pinch had computed it for $k \leq 35$.
\end{abstract}

\maketitle

\section{Introduction}\label{sec:intro}

\textit{Carmichael numbers} are the composite numbers $n$ satisfying $b^n \equiv b \pmod{n}$ for every integer $b$.
Fermat's little theorem says that the congruence must hold whenever $n$ is a prime, and so these numbers might also be called \textit{absolute Fermat pseudoprimes}.
In \cite{Carmichael1910}, Robert Carmichael proved that such numbers must be square-free and satisfy $(p-1) \mid (n-1)$ for every prime $p \mid n$, and he gave examples\footnote{
In 1899, Korselt had previously proven this result but gave no examples \cite{Korselt1899}.
In 1885, \v{S}imerka found the first seven examples but did not state how they were found \cite{Simerka1885}.}.

\begin{theorem}[Korselt's Criterion]
A composite number $n$ is a Carmichael number if and only if $n$ is square-free and $(p-1) \mid (n-1)$ for all prime divisors $p$ of $n$.
\end{theorem}

Every Carmichael number has at least three prime factors, and Alford, Granville, and Pomerance proved that there are infinitely many \cite{AGP94}, making rigorous a heuristic construction of Erd\H{o}s \cite{Erdos56}.
It is natural to ask for the smallest examples with a prescribed number of prime factors.
Let $S_k$ denote the least Carmichael number with exactly $k$ prime factors.
This is sequence A006931 in the OEIS \cite{OEIS}.
The present work extends the known values of $S_k$ and analyzes the cost of the algorithm that finds them.

\subsection{History}

The first few values of $S_k$ fall out of complete tabulations of Carmichael numbers.
\v{S}imerka's list contains $S_3 = 561$.
The values $S_4 = 41041$ and $S_5 = 825265$ appear in Poulet's 1938 table of base-$2$ pseudoprimes up to $10^8$, where the absolute pseudoprimes are marked \cite{Poulet1938}\footnote{
Poulet presented a table to $5 \cdot 10^7$ at the 1926 meeting of the French Association for the Advancement of Science in Lyon.
He states that it was published in two 1928 issues of \textit{Sphinx-Oedipe}, but he was himself unsure whether those issues were ever distributed.
We cite the 1938 paper as the first reliably accessible record, while noting that the computation is older.}.
In 1975, J.D. Swift found all Carmichael numbers up to $10^9$ \cite{Swift1975}, and this gives $S_6 = 321197185$.
Later tabulations refined Swift's prime-by-prime approach.
The tabulations of Pomerance, Selfridge, and Wagstaff to $25 \cdot 10^9$ \cite{PSW80} and of Jaeschke to $10^{12}$ \cite{Jaeschke1990} give $S_7$ and $S_8$, respectively.
Keller's tabulation to $10^{13}$ gives $S_9$ \cite{Keller1988}, and Guthmann later extended the tabulation to $10^{14}$ \cite{Guthmann1992}.
Pinch's tabulation to $10^{15}$ \cite{Pinch93}, together with his later extensions \cite{Pinch17,Pinch18}, most recently to $10^{21}$ \cite{PinchTable}, gives $S_{10}$ through $S_{12}$.

Beyond this point the problem changes character.
A complete tabulation up to $S_k$ is out of reach for even moderate $k$; $S_{20}$ already has $36$ digits.
Instead, one searches only for Carmichael numbers with exactly $k$ prime factors below some bound $B$, and lowers $B$ whenever a new example is found.
Pinch used this approach to compute $S_k$ for $k \leq 35$, proving each value minimal, and maintains a table of the results \cite{PinchTable}.
This paper is the full account of a project at Butler University that began in 2025.
Aidan Johnson and the second author, working as undergraduates, found $S_k$ for $36 \leq k \leq 39$ with an early version of the algorithm described here, extending Pinch's list by four entries.
Those values were presented as a poster \cite{MRC25} and announced in \cite{SW25}, so that the work in progress would be public.
We now describe the algorithm in full, analyze it, and extend the computation to $k = 64$.

\subsection{Contributions}

Our contributions are as follows.
We describe an algorithm for finding $S_k$ that refines the prime-by-prime approach.
Rather than appending primes until only one remains to be found, the search stops appending primes as soon as the preproduct $P$ satisfies $P\lambda(P) > B$, much as Pinch's later searches did \cite{Pinch17}.
The new ingredient is what happens next.
All of the remaining prime factors are recovered at once from the single candidate $n = PR$, with $R \equiv P^{-1} \pmod{\lambda(P)}$, using the polynomial-time algorithm of \cite{SW25}.
We prove an unconditional bound on the number of candidates this algorithm considers, in terms of how far $B$ exceeds the product of the first $k$ primes.
Under a natural conjecture on $S_k$, this gives a running time of $e^{O(k)} = S_k^{O(1/\log k)}$.
We show that the constant in the exponent is sharp for the square-free numbers the bound counts, so that any further savings must come from admissibility.
We also give a heuristic that makes a prediction we find especially pleasing: if the number of primes recovered by factoring grows sublinearly in $k$, as our data suggest, then the prime-by-prime search and ours agree to leading order, and the advantage of ours is subexponential.
Finally, we implemented the algorithm with a work-stealing parallel scheduler and used it to find $S_k$ for $40 \leq k \leq 64$.

The remainder of the paper is organized as follows.
In Section \ref{sec:algorithm}, we describe the search tree and the algorithm, and we compare it with the naive and prime-by-prime searches on the same tree.
In Section \ref{sec:analysis}, we analyze the algorithm as $k$ grows.
In Section \ref{sec:implementation}, we discuss the implementation.
In Section \ref{sec:results}, we give the new values of $S_k$ and statistics on the computation.
We conclude in Section \ref{sec:lpf} by comparing $S_k$ with the least cyclic number with $k$ prime factors, and by studying the least prime factor of $S_k$.
Along the way we answer a question of Nielsen about the structure of these cyclic numbers.

\section{The algorithm}\label{sec:algorithm}

Throughout, $p_i$ denotes the $i$th prime, $\theta(x)=\sum_{p\le x}\log p$,
and $\omega(n)$ is the number of distinct prime factors of $n$. For a
squarefree $P>1$ we write $\pmx(P)$ for its largest prime factor and
\[
  \lambda(P)=\lcm\{p-1 : p\mid P\}
\]
for Carmichael's function. By Korselt's criterion, a squarefree $n=PR$ with
$\gcd(P,R)=1$ can be a Carmichael number only if $\lambda(P)\mid n-1$, that
is,
\[
  R\equiv P^{-1}\pmod{\lambda(P)}.
\]
This requires $\gcd(P,\lambda(P))=1$.

A positive integer $n$ is \emph{cyclic} if $\gcd(n, \varphi(n)) = 1$.
(These are the $n$ for which every group of order $n$ is cyclic.)
A cyclic number is square-free, and for square-free $P$ the numbers $\lambda(P)$ and $\varphi(P)$ have the same prime factors.
So the condition above says that $P$ is cyclic, and in particular every Carmichael number is cyclic.
For a cyclic $P$ and a prime $q > \pmx(P)$, we say that $q$ is \emph{admissible} to $P$ if $Pq$ is cyclic, that is, if $q \not\equiv 1 \pmod p$ for every prime $p \mid P$.
Only admissible primes can be appended to $P$.

\subsection{The search tree}

\begin{definition}\label{def:tree}
Fix $B$ and $k$. The \emph{search tree} $\Tree(B,k)$ has root $P=1$. A node
$P$ with $\omega(P)=j<k$ has as children the nodes $Pq$, where $q$ runs over
the primes admissible to $P$ in the interval
\[
  \pmx(P) < q \le (B/P)^{1/(k-j)}
\]
(with $\pmx(1):=2$). We say $P$ has \emph{depth} $j=\omega(P)$.
\end{definition}

The upper endpoint guarantees that $P$ can still be completed by $k-j$
larger primes without exceeding $B$. Every cyclic $n\le B$ with
$\omega(n)=k$ appears exactly once as a depth-$k$
node, through its factors in increasing order.

\subsection{Three algorithms on one tree}

We compare three ways of searching $\Tree(B,k)$. They differ only in where
they stop enumerating primes explicitly.

\begin{itemize}
\item[$\mathsf{N}$] (naive) Explores all of $\Tree(B,k)$. Every depth-$k$ node
  is a complete candidate and is tested with Korselt's criterion.
\item[$\mathsf{P}$] (prime by prime) Explores depths $\le k-1$. At each
  depth-$(k-1)$ node $P$ it runs through $R\equiv P^{-1}\pmod{\lambda(P)}$ with
  $\pmx(P)<R\le B/P$, and tests $R$ for primality.
\item[$\mathsf{W}$] (this paper) A node $P$ of depth $j$ becomes a \emph{leaf} if
  $k-j\le 2$ or $P\lambda(P)>B$. At a leaf it runs through
  $R\equiv P^{-1}\pmod{\lambda(P)}$ with $\pmx(P)<R\le B/P$. For each such $R$ it
  asks whether $PR$ is a Carmichael number.
  Pseudocode is given in Algorithm~\ref{alg:W}.
\end{itemize}

\begin{algorithm}
\caption{Carmichael numbers $n \le B$ with $\omega(n) = k$}\label{alg:W}
\begin{algorithmic}[1]
\Require The tabulation bound $B$, a preproduct $P$ with largest prime $p$, $\lambda(P)$, and the number $m$ of primes still to be found.
\Function{Search}{$B, P, \lambda(P), m$}
  \If{$m \le 2$ \textbf{or} $P\lambda(P) > B$}
    \State $r^\star \gets P^{-1} \bmod \lambda(P)$
    \State $R \gets$ the least $R \equiv r^\star \pmod{\lambda(P)}$ with $R > p$
    \While{$PR \le B$}
      \If{$PR$ is a Carmichael number with $k$ prime factors}
        \State report $PR$ and \textbf{break}
      \EndIf
      \State $R \gets R + \lambda(P)$
    \EndWhile
  \Else
    \For{each prime $q$ admissible to $P$ in $\bigl(p, \sqrt[m]{B/P}\,\bigr]$}
      \State \Call{Search}{$B, Pq, \lcm(\lambda(P), q-1), m-1$}
    \EndFor
  \EndIf
\EndFunction
\end{algorithmic}
\end{algorithm}

The search begins with \textsc{Search}$(B, 1, 1, k)$, where we take $p = 2$ for $P = 1$.

The heart of the method is the test $P\lambda(P) > B$ in line~2.
Once it holds, every completion satisfies $R \leq B/P < \lambda(P)$, so $R = r^\star$ is the only possibility, and testing it directly is cheaper than branching further.
(The same observation, that at most $1 + x/(d\lambda(d))$ Carmichael numbers $n \leq x$ are divisible by $d$, appears in Pomerance's account of the upper bound for the number of Carmichael numbers \cite[(5.3)]{Pomerance89}.)
Some branches never reach $P\lambda(P) > B$, and the test $m \leq 2$ ensures that these end in a leaf too.
At every leaf, the loop in lines 5--10 runs over all completions of $P$ up to $B$; when $P\lambda(P) > B$ it runs at most once.
So the test in line~2 affects only the cost of the search and not its correctness, which gives us freedom in how we implement it (see Section~\ref{sec:implementation}).
Algorithm $\mathsf P$ is obtained by replacing the test in line~2 with $m = 1$, in which case line~6 only needs to ask whether $R$ is prime.

Pinch's later searches \cite{Pinch17,Pinch18} go further than $\mathsf P$.
They stop at a level where $\lambda(P)$ is large enough to limit the possible cofactors, and factor those cofactors directly, so they visit essentially the tree of $\mathsf W$.
The candidate counts of Section~\ref{sec:analysis} apply to them as well.
Where they differ is line~6.
Pinch did not specify a factoring method, and for the purposes of analysis any method can be used, but no general factoring algorithm is known to run in polynomial time.
In $\mathsf W$ the cofactor is not factored in isolation.
Line~6 works with $n = PR$ itself, and the algorithm of \cite{SW25} answers it in polynomial time.
So for $\mathsf W$ the bounds of Section~\ref{sec:analysis} bound the running time, and not only the number of candidates.

Line~6 is an instance of the decision problem \textsc{carmichaels}: given $n$, is $n$ a Carmichael number?
This is answered in polynomial time by \cite{SW25}, and we discuss how in Section~\ref{sec:factoring}.
If the answer is yes, then $n - 1$ is a multiple of $\lambda(n)$, which lets us factor $n$ completely, and so counting its prime factors costs nothing extra.
Since $R$ increases along the progression, the first Carmichael number found at a leaf is the least one there, and there is no reason to continue.

\begin{definition}
At a leaf $P$ of depth $j$, the factoring step recovers the remaining $k - j$ primes.
Let $t_k$ be the average of $k - j$ over the leaves of $\mathsf W$ with $P\lambda(P) > B$, in a search with $B$ close to $S_k$.
\end{definition}

Algorithm $\mathsf P$ corresponds to always recovering exactly one prime.
The gain of $\mathsf W$ comes from $t_k$ being larger.
Section~\ref{sec:analysis} shows $t_k < k/2$, and in our computations $t_k$ grows slowly with $k$ (Table~\ref{tab:tk}).

\subsection{The Carmichael query}\label{sec:factoring}

We answer line~6 with the algorithm of \cite{SW25}, which uses Fermat tests modulo $n$ and GCDs to factor $R$.
Its correctness is unconditional.
It runs in deterministic polynomial time under ERH, and in expected polynomial time unconditionally if the bases are chosen at random \cite[Thm.~10]{SW25}.
We describe our implementation of it in Section~\ref{sec:implementation}.

\section{Analysis of the algorithm}\label{sec:analysis}

Every node visited by $\mathsf N$, $\mathsf P$, or $\mathsf W$ is a node of
$\Tree(B,k)$. The following two constraints on nodes are what drive
everything below.

\begin{lemma}\label{lem:node}
Let $P\ne1$ be a node of $\Tree(B,k)$ of depth $j$. Then
\begin{enumerate}
\item[(i)] $P\,\pmx(P)^{\,k-j}\le B$;
\item[(ii)] $P\le B^{j/k}$.
\end{enumerate}
\end{lemma}

\begin{proof}
(i) $P=P'q$ with $q=\pmx(P)$ a child of $P'$, which has depth $j-1$, so
$q^{k-j+1}\le B/P'$; that is, $P q^{k-j}\le B$.
(ii) Since $P\le \pmx(P)^j$, we have
$P^{(k-j)/j}\le\pmx(P)^{k-j}$, so by (i)
$P^{k/j}=P\cdot P^{(k-j)/j}\le P\,\pmx(P)^{k-j}\le B$.
\end{proof}

\subsection{An unconditional bound}

We want to count the nodes of $\Tree(B,k)$.
By Lemma~\ref{lem:node}(ii), a node of depth $j$ is a square-free number up to $B^{j/k}$ with exactly $j$ prime factors.
So it is enough to count those numbers, and here we ignore admissibility entirely.
The argument has three steps.
First, Rankin's trick (Lemma~\ref{lem:rankin}) bounds the count by an Euler product $\prod_p(1+zp^{-\sigma})$, with two free parameters $\sigma$ and $z$.
Second, we estimate the logarithm of that product.
It is a sum over primes, and the primes that matter lie near $p_k$, where there are about $k$ of them per unit of $p/p_k$.
So the sum is $k$ times an integral, up to $o(k)$ (Lemma~\ref{lem:sums}(b)).
Third, we choose $z$ and $\sigma$ to make the bound as small as possible.
What comes out is $e^{(c(\beta)+o(1))k}$, where $\beta$ measures how far $B$ exceeds the product of the first $k$ primes, and $c(\beta)$ is defined by that optimization (Theorem~\ref{thm:main}).

\begin{lemma}[Rankin's bound \cite{Rankin38}]\label{lem:rankin}
Let $\mathcal{Q}$ be a set of primes, $x\ge1$, $j\ge0$, and let
$N_j(x;\mathcal{Q})$ be the number of squarefree $n\le x$ with $\omega(n)=j$
and all prime factors in $\mathcal{Q}$. For every $\sigma\ge0$ and $z>0$,
\[
  N_j(x;\mathcal{Q})\;\le\;\sum_{\substack{n\ \text{counted}}}\Big(\frac{x}{n}\Big)^{\sigma}
  \;\le\; x^{\sigma}\,z^{-j}\prod_{p\in\mathcal{Q}}\bigl(1+z\,p^{-\sigma}\bigr).
\]
\end{lemma}

\begin{proof}
Each counted $n$ has $(x/n)^\sigma\ge1$. The middle sum is at most $x^\sigma$
times the coefficient of $z^j$ in $\prod_{p\in\mathcal Q}(1+zp^{-\sigma})$,
and since every coefficient is nonnegative, that coefficient is at most
$z^{-j}$ times the whole product.
\end{proof}

This is Rankin's method.
Pomerance \cite{Pomerance89} gives an exposition with many applications, among them an upper bound for the number of Carmichael numbers.

Write
\begin{equation}\label{eq:beta}
  \log B=\theta(p_k)+\beta k .
\end{equation}
Since $\theta(p_k)=k\log k+k\log\log k-k+o(k)$ is the logarithm of the
product of the first $k$ primes, $\beta$ measures the \emph{excess} of the
search bound over the smallest conceivable $k$-factor product, per prime.
For $\sigma>1$ put
\[
  A_\sigma=\frac{\pi}{\sin(\pi/\sigma)},\qquad
  g(\sigma)=\frac{\sigma}{A_\sigma}=\frac{\sigma}{\pi}\sin\frac{\pi}{\sigma}\in(0,1).
\]
The constant $A_\sigma$ is the integral $\int_0^\infty \log(1+u^{-\sigma})\,du$ from the second step, and we evaluate it in the proof of Lemma~\ref{lem:sums}(b).
The closed form is not needed for the bound itself.
But it makes $c(\beta)$ explicit, and its expansion for large $\sigma$ is the source of the Hardy--Ramanujan constant in the remark on small excess below.
Define, for $\beta\ge0$,
\begin{equation}\label{eq:c}
  c(\beta)=\inf_{\sigma>1}\ \sigma\bigl(\beta-\log g(\sigma)\bigr).
\end{equation}

\begin{lemma}\label{lem:sums}
As $k\to\infty$:
\begin{enumerate}
\item[(a)] $k\log p_k-\theta(p_k)=(1+o(1))\,k$;
\item[(b)] for fixed $\sigma>1$,
  $\displaystyle\sum_{p}\log\Bigl(1+w\Bigl(\frac{p_k}{p}\Bigr)^{\sigma}\Bigr)
  =\bigl(A_\sigma w^{1/\sigma}+o(1)\bigr)k$, the sum over all primes, uniformly for $0<w\le1$.
\end{enumerate}
\end{lemma}

\begin{proof}
(a) By partial summation,
$\theta(x)-\pi(x)\log x=-\int_2^x\pi(t)\,dt/t=-(1+o(1))x/\log x$. Take
$x=p_k$ and use $p_k/\log p_k\sim k$.

(b) Put $h(u)=\log(1+wu^{-\sigma})$, so the sum is
$\sum_p h(p/p_k)$. Fix $\delta\in(0,1)$. Primes $p\le p_k^{1-\delta}$
contribute $O\big(\pi(p_k^{1-\delta})\log p_k\big)=o(k)$, since
$h(u)\le\log(1+w)+\sigma\log(1/u)$ for $u\le1$. Primes $p>p_k^{1+\delta}$
contribute at most $w\,p_k^\sigma\sum_{p>p_k^{1+\delta}}p^{-\sigma}
=O\big(p_k^{1-\delta(\sigma-1)}\big)=o(k)$. On the middle range,
$d\pi(t)=(1+O(\delta))\,dt/\log p_k$ by the prime number theorem, so the
contribution is $(1+O(\delta))\frac{p_k}{\log p_k}\int_0^\infty h(u)\,du$.
Integrating by parts,
$\int_0^\infty\log(1+u^{-\sigma})\,du=\sigma\int_0^\infty\frac{du}{1+u^\sigma}
=\frac{\pi}{\sin(\pi/\sigma)}$, where the last step substitutes $t=u^\sigma$ and uses the Beta integral $\int_0^\infty t^{s-1}(1+t)^{-1}\,dt=\Gamma(s)\Gamma(1-s)$ with $s=1/\sigma$, together with Euler's reflection formula. Scaling gives
$\int_0^\infty h=A_\sigma w^{1/\sigma}$. Let $\delta\to0$.
Each error term above is bounded independently of $w\in(0,1]$, so the estimate is uniform in $w$.
\end{proof}

\begin{theorem}\label{thm:main}
Fix $\beta\ge0$ and suppose $\log B\le\theta(p_k)+\beta k$. Then, as
$k\to\infty$:
\begin{enumerate}
\item[(i)] the number of depth-$k$ nodes of $\Tree(B,k)$ is at most
  $e^{(c(\beta)+o(1))k}$;
\item[(ii)] $|\Tree(B,k)|\le e^{(c(\beta+1/e)+o(1))k}$;
\item[(iii)] each of $\mathsf N$, $\mathsf P$, $\mathsf W$ examines at most
  $e^{(c(\beta+1/e)+o(1))k}$ candidates (tree nodes plus arithmetic-progression
  terms), and every prime it enumerates is at most $p_k e^{\beta k}$.
\end{enumerate}
\end{theorem}

\begin{proof}
Fix $\sigma>1$, let $0<w\le1$, and apply Lemma~\ref{lem:rankin} with
$\mathcal{Q}$ the set of all primes and $z=w\,p_k^{\sigma}$.

\smallskip\noindent
\emph{Depth $j=\alpha k$.} By Lemma~\ref{lem:node}(ii) the depth-$j$ nodes are
counted by $N_j(B^{j/k})$, and
\[
  \log N_j(B^{j/k})\le \sigma j\Big(\frac{\log B}{k}-\log p_k\Big)-j\log w
  +\sum_p\log\Bigl(1+w\Bigl(\frac{p_k}{p}\Bigr)^\sigma\Bigr).
\]
By hypothesis and Lemma~\ref{lem:sums}(a),
$\log B/k-\log p_k\le\beta-1+o(1)$. By Lemma~\ref{lem:sums}(b), the right-hand side
is at most $k\big[\alpha\sigma(\beta-1)-\alpha\log w+A_\sigma w^{1/\sigma}\big]+o(k)$.
Choosing $w^{1/\sigma}=\alpha g(\sigma)$ (so $w\le1$, since $\alpha\le1$ and $g(\sigma)<1$; the $o(1)$ is uniform in $\alpha$ by Lemma~\ref{lem:sums}(b)) and using $A_\sigma g(\sigma)=\sigma$
gives
\[
  \log N_j(B^{j/k})\le k\,\alpha\sigma\bigl(\beta+\log(1/\alpha)-\log g(\sigma)\bigr)+o(k).
\]
For $\alpha=1$ this is $k\sigma(\beta-\log g(\sigma))+o(k)$. For general
$\alpha\in(0,1]$ use $\alpha\log(1/\alpha)\le1/e$. Summing over the $k+1$
depths only adds $O(\log k)$. Since $\sigma>1$ was arbitrary, (i) and (ii)
follow.

\smallskip\noindent
\emph{Progression terms.} For $\mathsf P$ and $\mathsf W$, a leaf $P$ yields at
most $1+B/(P\lambda(P))$ terms (compare \cite[(5.3)]{Pomerance89}). If $P\lambda(P)>B$ this is at most $2$, and
those leaves are already counted in (ii). Otherwise $P$ has depth $j\ge k-2$.
Then $\lambda(P)\ge\pmx(P)-1\ge\pmx(P)/2$, and $\rho:=B/(P\,\pmx(P))\ge1$ by
Lemma~\ref{lem:node}(i). So $\rho\le\rho^\sigma\le(B/(Pp_j))^\sigma$,
because $\pmx(P)\ge p_j$. Summing over all squarefree $P$ with $\omega(P)=j$ and
applying the middle inequality of Lemma~\ref{lem:rankin} gives
\begin{multline*}
  \log\sum_{P}\frac{B}{P\,\pmx(P)}
  \le\sigma(\log B-k\log p_k)+\sigma(k-j)\log p_k-\sigma\log p_j \\
  -j\log w+\bigl(A_\sigma w^{1/\sigma}+o(1)\bigr)k .
\end{multline*}
Since $k-j\le2$, this is the $\alpha=1$ case of the bound above up to
$O(\log k)$.

\smallskip\noindent
\emph{Largest prime.} A child $q$ of a depth-$(j-1)$ node $P'$ satisfies
$q^{k-j+1}\le B/P'\le B/(p_1\cdots p_{j-1})$. Hence
$(k-j+1)\log q\le\theta(p_k)-\theta(p_{j-1})+\beta k\le(k-j+1)\log p_k+\beta k$,
so $\log q\le\log p_k+\beta k$.
\end{proof}

\begin{remark}\label{rem:binomial}
Taking $\sigma=0$ in Lemma~\ref{lem:rankin} recovers, up to a factor
$|\mathcal Q|+1$, the binomial count $\binom{|\mathcal Q|}{j}$ of $j$-subsets
of a finite set of primes $\mathcal Q$. That count ignores the constraint
$n\le B$, which is what actually limits the search. This is why the proof of
Theorem~\ref{thm:main} uses a constant $\sigma>1$.
\end{remark}

\begin{remark}
Theorem~\ref{thm:main} uses neither admissibility nor any information about
$\lambda$ beyond the trivial $\lambda(P)\ge\pmx(P)-1$. It holds for any
search of this shape. Admissibility and the actual size of $\lambda(P)$ can
only improve the constants.
\end{remark}

\subsection{Finding $S_k$ without knowing it}

Theorem~\ref{thm:main} bounds the cost of searching $\Tree(B,k)$ for any $B$.
It does not assume that $B \geq S_k$.
A search with $B < S_k$ simply finds nothing, and a search with $B \geq S_k$ finds $S_k$.
To find $S_k$ we must therefore choose $B$, and we do not know $S_k$ in advance.
One might worry that the analysis is circular: we use a heuristic to choose $B$, and then analyze the search with that $B$.
It is not, because the procedure we actually run does not need to know $S_k$.

\begin{definition}[Doubling passes]\label{def:doubling}
Let $B_0$ be the product of the first $k$ odd primes.
For $i = 0, 1, 2, \ldots$, search $\Tree(2^i B_0, k)$, and stop at the first pass that finds a Carmichael number with $k$ prime factors.
Report the least one found in that pass.
\end{definition}

Every Carmichael number is odd and square-free, so $S_k \geq B_0$.
Pass $i$ finds $S_k$ exactly when $2^i B_0 \geq S_k$, and so the procedure stops at the least such $i$ and reports $S_k$.
In particular the final bound satisfies $B < 2S_k$.
(The factor $2$ is not special; with any constant ratio $r > 1$ the final bound satisfies $B < rS_k$, which changes $\beta$ by only $(\log r)/k$.)
If no Carmichael number with $k$ prime factors exists, then the procedure does not stop.
This is not a concern in our range: Alford, Grantham, Hayman, and Shallue constructed Carmichael numbers with exactly $k$ prime factors for every $3 \leq k \leq 19\,565\,220$ \cite{AGHS14}.

The cost of the procedure is governed by how far $S_k$ lies above the product of the first $k$ primes.
Define $\beta_k$ by
\[
  \log S_k = \theta(p_k) + \beta_k k ,
\]
so that $\beta_k$ is the quantity $\beta$ of \eqref{eq:beta} evaluated at $B = S_k$.
This is a property of the sequence $S_k$ alone, and it involves no choice on our part.

\begin{corollary}\label{cor:runtime}
Let $\beta^* = \limsup_{k \to \infty} \beta_k$, and suppose $\beta^* < \infty$.
Then running any of $\mathsf N$, $\mathsf P$, $\mathsf W$ in doubling passes finds $S_k$ after examining at most
\[
  e^{(c(\beta^* + 1/e) + o(1))k}
\]
candidates in total.
Since $\log S_k \geq (1 + o(1))k \log k$, the running time is
\[
  T(k) = e^{O(k)} \cdot \mathrm{poly}(\log S_k) = S_k^{\,O(1/\log k)} .
\]
Here the per-candidate cost is that of a Korselt check for $\mathsf N$, a deterministic primality test for $\mathsf P$ (unconditionally polynomial, by \cite{AKS}), and the factoring step of \cite{SW25} for $\mathsf W$.
For $\mathsf W$ this bound is unconditional for the expected running time of the randomized factoring step; see Section~\ref{sec:factoring}.
\end{corollary}

\begin{proof}
Fix $\varepsilon > 0$.
For all large $k$ we have $\beta_k \leq \beta^* + \varepsilon$.
Every pass has $B < 2S_k$, so
\[
  \log B < \theta(p_k) + \beta_k k + \log 2 \leq \theta(p_k) + (\beta^* + 2\varepsilon)k
\]
for large $k$.
By Theorem~\ref{thm:main}(iii), each pass examines at most $e^{(c(\beta^* + 2\varepsilon + 1/e) + o(1))k}$ candidates.
The number of passes is at most $1 + \log_2(S_k/B_0)$.
Since $\log B_0 = \theta(p_{k+1}) - \log 2 \geq \theta(p_k) - \log 2$, this is at most $1 + ((\beta^* + \varepsilon)k + \log 2)/\log 2 = O(k)$.
The factor $O(k)$ is absorbed into the $o(1)$ in the exponent.
Finally, $c$ is continuous, so we may let $\varepsilon \to 0$.
The per-candidate cost is polynomial in $\log S_k = O(k \log k)$, which is $e^{O(\log k)}$.
\end{proof}

The bound in Corollary~\ref{cor:runtime} is \textit{output-sensitive}.
It is stated in terms of the answer $S_k$, through $\beta^*$, but the procedure itself uses no information about $S_k$.
The only hypothesis is that $\beta^*$ is finite, which is a statement about the size of $S_k$ and not about the algorithm.
The bound also covers the table of primes the search needs.
By the proof of Theorem~\ref{thm:main}, no node contains a prime larger than $p_k e^{\beta k}$, and a table of the primes up to that size costs $e^{(\beta + o(1))k}$.
Since $c(\beta) > \beta$, this does not change the exponent.

In practice we start at a bound $B_0$ given by a heuristic formula, rather than at the product of the first $k$ odd primes.
This changes only which of the cheap early passes are skipped.
If the formula undershoots $S_k$, nothing in the analysis changes.
If it overshoots, say $\log B_0 = \log S_k + \delta k$, then the single pass costs $e^{(c(\beta_k + \delta + 1/e) + o(1))k}$.
So a heuristic that overshoots by a factor of $e^{o(k)}$ does not change the exponent.

It remains to say why $\beta^*$ should be finite.
For $k \geq 3$, let $q$ be the least number such that the interval $(q/2, q]$ contains $k$ primes, and let $U(k)$ be the product of those $k$ primes.

\begin{conjecture}\label{conj:U}
$S_k<U(k)$ for all $k\ge3$.
\end{conjecture}

\begin{proposition}
$\log U(k)\le\theta(p_k)+(1+\log2+o(1))k$.
Consequently, Conjecture~\ref{conj:U} implies $\beta^* \leq 1 + \log 2$, and so $T(k) = e^{O(k)}$.
\end{proposition}

\begin{proof}
By the prime number theorem $\pi(q)-\pi(q/2)\sim q/(2\log q)$, so this
$q$ satisfies $q\sim2k\log k$. Every factor of $U(k)$ is at most
$q$, so $\log U(k)\le k\log q=k\log k+k\log\log k+k\log2+o(k)$. Compare with
$\theta(p_k)=k\log k+k\log\log k-k+o(k)$.
\end{proof}

\subsection{The sharp constant}\label{sec:sharp}

Theorem~\ref{thm:main} is only an upper bound, and Rankin's method can be generous.
Here we show that it is not, for the set the theorem counts: $c(\beta)$ is the true constant for the square-free numbers with $k$ prime factors.
So the gap between the theorem and our data (Section~\ref{sec:practice}) comes from admissibility, which the theorem ignores, and not from Rankin's method.

The optimal $\sigma$ and $z$ in the proof of Theorem~\ref{thm:main} form a saddle point.
For smooth numbers \cite{HT86} and for integers with a given number of prime factors \cite{HT88} (see also \cite{Tenenbaum}), the saddle-point method shows that an optimized Rankin bound has the correct exponential order.
Our regime lies beyond the range of those results: the formula of \cite{HT88} for $\#\{n \leq x : \omega(n) = \nu\}$ holds for $\nu \leq \log x/(\log\log x)^2$, while we need $\nu = k \sim \log x/\log\log x$.
But at exponential scale a direct construction suffices.
We choose a fixed fraction of the primes in each of a bounded number of ranges near $p_k$, and count the choices with binomial coefficients.
The fractions come from the optimal parameters in Rankin's bound.
Because the ranges have length proportional to $p_k$, the prime number theorem is all we need.
Erd\H{o}s and Nicolas \cite{EN81} used the same kind of construction for the case $\beta k = \log p_{k+1}$.

Write $H(t) = -t\log t - (1-t)\log(1-t)$ for $0 \leq t \leq 1$.

\begin{theorem}\label{thm:sharp}
Fix $\beta > 0$, and let $M_k(\beta)$ be the number of square-free $n \leq e^{\theta(p_k)+\beta k}$ with $\omega(n) = k$.
Then $\log M_k(\beta) = \bigl(c(\beta) + o(1)\bigr)k$ as $k \to \infty$.
\end{theorem}

The upper bound is the case $\alpha = 1$ of the proof of Theorem~\ref{thm:main}, which counts all square-free $n$ with $\omega(n) = k$ and never uses admissibility.
The lower bound rests on two lemmas.

\begin{lemma}\label{lem:profile}
Fix $\beta > 0$.
The infimum defining $c(\beta)$ is attained at some $\sigma_0 > 1$.
Put $w_0 = g(\sigma_0)^{\sigma_0}$ and $f(u) = 1/(1 + u^{\sigma_0}/w_0)$ for $u > 0$.
Then
\begin{gather*}
  \text{(a)}\ \int_0^\infty f(u)\,du = 1, \qquad
  \text{(b)}\ \int_0^\infty f(u)\log u\,du = \beta - 1, \\
  \text{(c)}\ \int_0^\infty H\bigl(f(u)\bigr)\,du = c(\beta).
\end{gather*}
\end{lemma}

\begin{proof}
Let $\Phi(\sigma) = \sigma(\beta - \log g(\sigma))$.
As $\sigma \to 1^+$ we have $g(\sigma) \to 0$, and $\Phi(\sigma) \geq \sigma\beta$ since $g < 1$, so $\Phi \to \infty$ at both ends of $(1,\infty)$.
It is continuous, so its infimum is attained at some $\sigma_0$.

For $\sigma > 1$ and $w > 0$ put $f_{\sigma,w}(u) = wu^{-\sigma}/(1 + wu^{-\sigma})$, so that $f = f_{\sigma_0,w_0}$, and put
\[
  I(\sigma,w) = \int_0^\infty \log(1 + wu^{-\sigma})\,du = A_\sigma w^{1/\sigma},
  \qquad
  F(\sigma,w) = \sigma(\beta - 1) - \log w + I(\sigma,w),
\]
where the closed form for $I$ is the computation in the proof of Lemma~\ref{lem:sums}(b).
For fixed $\sigma$, $F$ is a convex function of $\log w$, minimized where $w^{1/\sigma} = \sigma/A_\sigma = g(\sigma)$, and there $F = \Phi(\sigma)$.
So $\Phi(\sigma) = \min_w F(\sigma, w)$, with minimizer $w^*(\sigma) = g(\sigma)^\sigma$, and $w^*(\sigma_0) = w_0$.

(a) $\int_0^\infty f_{\sigma,w}\,du = \partial I/\partial\log w = w^{1/\sigma}A_\sigma/\sigma$, which is $1$ at $w = w^*(\sigma)$.

(b) Differentiating under the integral sign (by dominated convergence, for $\sigma$ in a compact subset of $(1,\infty)$), $\partial I/\partial\sigma = -\int_0^\infty f_{\sigma,w}(u)\log u\,du$.
Since $\partial F/\partial w = 0$ at $w = w^*(\sigma)$, we get $\Phi'(\sigma) = \beta - 1 - \int_0^\infty f_{\sigma,w^*(\sigma)}(u)\log u\,du$.
At the minimum $\Phi'(\sigma_0) = 0$, which is (b).

(c) Write $a = w_0u^{-\sigma_0}$, so $f = a/(1+a)$, $\log f = \log a - \log(1+a)$, and $\log(1-f) = -\log(1+a)$.
Then
\[
  H(f) = \log(1+a) - f\log a = \log(1 + w_0u^{-\sigma_0}) + f(u)\bigl(\sigma_0\log u - \log w_0\bigr).
\]
Integrating and using (a) and (b), $\int_0^\infty H(f) = I(\sigma_0,w_0) + \sigma_0(\beta - 1) - \log w_0 = F(\sigma_0,w_0) = c(\beta)$.
\end{proof}

The integrals in Lemma~\ref{lem:profile} converge absolutely.
Near $0$ we have $f \leq 1$ and $H \leq \log 2$, and for large $u$ we have $f(u) \leq w_0u^{-\sigma_0}$ and $H(f(u)) = O(u^{-\sigma_0}\log u)$, with $\sigma_0 > 1$.

\begin{lemma}\label{lem:construct}
Fix $\beta > 0$ and $\eta > 0$.
For all large $k$ there are at least $e^{(c(\beta) - \eta)k}$ square-free $n$ with $\omega(n) = k$ and $\bigl|\log n - \theta(p_k) - \beta k\bigr| \leq \eta k$.
\end{lemma}

\begin{proof}
Let $f$ be as in Lemma~\ref{lem:profile}.
Choose parameters $0 < \varepsilon < 1 < 2 < U$ and $\delta > 0$, in that order, as described at the end of the proof.
Let $I_1, \dots, I_L$ be the consecutive intervals of length $\delta$ that partition $[\varepsilon, U)$, taking $(U - \varepsilon)/\delta$ to be an integer, and let $u_i$ be the left endpoint of $I_i$.
Every $n$ we build is the product of all primes $p \leq \varepsilon p_k$ and, for each $i$, a set of $m_i$ primes $p$ with $p/p_k \in I_i$.
These primes are distinct, so $n$ is square-free.
Let $N_i$ be the number of primes $p$ with $p/p_k \in I_i$, and let $s = \pi(\varepsilon p_k)$.
Since the $I_i$ lie in the fixed range $[\varepsilon, U]$, the prime number theorem and $p_k/\log p_k \sim k$ give
$N_i = \delta k\,(1 + o(1))$, uniformly in $i$, and $s = \varepsilon k\,(1 + o(1))$.

\emph{Choice of the $m_i$.}
Start from $m_i^0 = \lfloor f(u_i)N_i \rfloor$.
Since $f$ is decreasing, $\delta\sum_i f(u_i)$ differs from $\int_\varepsilon^U f$ by at most $\delta$, so by Lemma~\ref{lem:profile}(a)
\[
  s + \sum_i m_i^0 = k\Bigl(1 + \varepsilon - \int_0^\varepsilon f - \int_U^\infty f + O(\delta)\Bigr) + o(k).
\]
So the defect $d = k - s - \sum_i m_i^0$ satisfies $|d| \leq \tau k + o(k)$, where $\tau = \varepsilon + \int_U^\infty f + \delta$.
Let $J$ be the set of $i$ with $I_i \subset [1,2]$, and let $\rho = \min(f(2), 1 - f(1)) > 0$.
For $i \in J$ the block has at least $\rho N_i - 1$ unused primes and at least $\rho N_i - 1$ chosen ones.
If $\tau < \rho/4$, we can absorb $d$ by changing each $m_i^0$ with $i \in J$ by at most $2\tau N_i + 1$, in the direction of the sign of $d$.
This gives $m_i$ with $0 \leq m_i \leq N_i$ and $s + \sum_i m_i = k$, so that every $n$ built has $\omega(n) = k$, and with
$|m_i/N_i - f(u_i)| \leq 1/N_i$ for $i \notin J$ and $|m_i/N_i - f(u_i)| \leq 2\tau + 2/N_i$ for $i \in J$.

\emph{Counting.}
The number of such $n$ is $\prod_i \binom{N_i}{m_i}$, and $\log\binom{N}{m} \geq N H(m/N) - \log(N+1)$.
Let $\mu$ be a modulus of continuity for $H$ on $[0,1]$.
Then
\begin{align*}
  \log \prod_i \binom{N_i}{m_i}
  &\geq \sum_i N_i H\bigl(f(u_i)\bigr) - k\,\mu(3\tau) - O(L\log k) \\
  &= k\Bigl(\delta\sum_i H\bigl(f(u_i)\bigr) - \mu(3\tau)\Bigr) + o(k).
\end{align*}
Here the term $k\,\mu(3\tau)$ comes from the blocks in $J$, which contain at most $(1 + o(1))k$ primes in all, and every other block loses only $N_i\,\mu(1/N_i) = O(\log k)$.
As $\delta \to 0$, $\delta\sum_i H(f(u_i)) \to \int_\varepsilon^U H(f)$, and by Lemma~\ref{lem:profile}(c) this is $c(\beta)$ minus the two tails $\int_0^\varepsilon H(f) \leq \varepsilon\log 2$ and $\int_U^\infty H(f)$.

\emph{Size.}
Since exactly $k$ primes divide $n$,
\[
  \log n - \theta(p_k) = \bigl(k\log p_k - \theta(p_k)\bigr) + \sum_{p \leq \varepsilon p_k}\log\frac{p}{p_k} + \sum_i \sum_{\substack{p \mid n \\ p/p_k \in I_i}} \log\frac{p}{p_k}.
\]
The first term is $k + o(k)$ by Lemma~\ref{lem:sums}(a).
For the second, the partial summation in the proof of Lemma~\ref{lem:sums}(a), applied at $x = \varepsilon p_k$, gives
\[
  \sum_{p \leq \varepsilon p_k}\log\frac{p}{p_k} = \theta(\varepsilon p_k) - \pi(\varepsilon p_k)\log p_k = -\varepsilon\bigl(1 + \log(1/\varepsilon)\bigr)k + o(k).
\]
For the third, each prime counted in $I_i$ has $|\log(p/p_k) - \log u_i| \leq \delta/\varepsilon$, and there are at most $k$ of them.
Also $|\log u_i| \leq \Lambda = \max(\log(1/\varepsilon), \log U)$ for all $i$, and $|\log u_i| \leq \log 2$ for $i \in J$.
So the third term is
\[
  \sum_i m_i\log u_i + O\Bigl(\frac{\delta}{\varepsilon}k\Bigr)
  = k\,\delta\sum_i f(u_i)\log u_i + O\Bigl(\tau k + L\Lambda + \frac{\delta}{\varepsilon}k\Bigr) + o(k),
\]
and $\delta\sum_i f(u_i)\log u_i \to \int_\varepsilon^U f\log u$ as $\delta \to 0$.
By Lemma~\ref{lem:profile}(b), $\int_\varepsilon^U f\log u = \beta - 1 - \int_0^\varepsilon f\log u - \int_U^\infty f\log u$, and both tails tend to $0$ as $\varepsilon \to 0$ and $U \to \infty$.
Putting the three terms together, $\log n - \theta(p_k) - \beta k$ is $k$ times a quantity that tends to $0$ as $\varepsilon \to 0$, then $U \to \infty$, then $\delta \to 0$, plus $o(k)$.

\emph{Choice of parameters.}
Given $\eta$, choose $\varepsilon$ small and then $U$ large, so that $\varepsilon + \int_U^\infty f < \rho/8$ and so that $\varepsilon(1 + \log(1/\varepsilon))$, the four tails, and $\mu(3\tau)$ are all small compared with $\eta$.
Then choose $\delta$ small compared with $\varepsilon\eta$ and with $\rho$.
For all large $k$, both the count and the size then satisfy the bounds claimed.
\end{proof}

\begin{proof}[Proof of Theorem~\ref{thm:sharp}]
It remains to prove the lower bound.
The function $c$ is an infimum of functions that are affine in $\beta$, so it is concave, and hence continuous, on $(0,\infty)$.
Given $\eta \in (0, \beta/2)$, apply Lemma~\ref{lem:construct} with $\beta - \eta$ in place of $\beta$ and $\eta/2$ in place of $\eta$.
Every $n$ it produces has $\log n \leq \theta(p_k) + (\beta - \eta/2)k$, so $\log M_k(\beta) \geq (c(\beta - \eta) - \eta/2)k$ for all large $k$.
Let $\eta \to 0$.
\end{proof}

\begin{remark}[Small excess]
As $\beta \to 0$, we have $c(\beta) \sim \pi\sqrt{2\beta/3}$, the Hardy--Ramanujan constant \cite{HR18}.
That is what one expects if the primes near $p_k$ were evenly spaced, since then the numbers counted in Theorem~\ref{thm:sharp} would correspond to partitions of about $\beta k^2$.
Erd\H{o}s and Nicolas \cite[p.~10]{EN81} use the same partition heuristic.
\end{remark}

\subsection{The number of primes recovered by factoring}

\begin{proposition}\label{prop:tk}
For $k>4$, every leaf of $\mathsf W$ has fewer than $k/2$ primes remaining. In particular $t_k<k/2$.
\end{proposition}

\begin{proof}
Let $P \neq 1$ be a leaf of depth $j$, so that $k - j$ primes remain.
If the leaf arose from the rule $k-j\le2$, then $k-j\le2<k/2$.
Otherwise $P\lambda(P)>B$.
By Lemma~\ref{lem:node}(i), and since $\lambda(P)\mid\varphi(P)$ gives $\lambda(P)<P$,
\[
  \pmx(P)^{k-j}\le B/P<\lambda(P)<P\le\pmx(P)^{j}.
\]
Hence $k-j<j$, that is, $k-j<k/2$.
\end{proof}

Table~\ref{tab:tk} shows $t_k$ as measured in our computations.
It grows slowly with $k$, from about $3$ at $k = 20$ to about $7$ at $k = 63$, and stays far below the bound $k/2$ of Proposition~\ref{prop:tk}.
Algorithm $\mathsf P$ recovers exactly one prime per leaf.
For our largest and most expensive searches, $\mathsf W$ recovers on average six to seven times as many.
We have no proof that $t_k$ keeps growing.

\begin{table}[ht]
\centering
\caption{The average number $t_k$ of primes recovered by factoring, over all leaves with $P\lambda(P) > B$, measured during clean passes.}\label{tab:tk}
\begin{tabular}{rr@{\qquad}rr@{\qquad}rr@{\qquad}rr}
\toprule
$k$ & $t_k$ & $k$ & $t_k$ & $k$ & $t_k$ & $k$ & $t_k$ \\
\midrule
20 & 2.924 & 40 & 5.726 & 46 & 5.984 & 54 & 6.562 \\
24 & 3.486 & 41 & 5.651 & 47 & 6.121 & 55 & 6.655 \\
29 & 4.145 & 42 & 5.805 & 48 & 6.253 & 56 & 6.752 \\
34 & 5.015 & 43 & 5.927 & 49 & 6.180 & 59 & 6.968 \\
36 & 5.460 & 44 & 5.707 & 50 & 6.275 & 60 & 6.990 \\
39 & 5.396 & 45 & 5.936 & 51 & 6.404 & 63 & 7.211 \\
\bottomrule
\end{tabular}
\end{table}

The bound $B$ in each run was close to, but not exactly, $S_k$, and the runs for different $k$ used different versions of our code.
This may account for the small departures from monotonicity (for example at $k = 36$ and $k = 44$).
\subsection{Where $\mathsf N$, $\mathsf P$, and $\mathsf W$ differ}

One might expect the three algorithms to differ the way $\binom{S_{\mathsf N}}{k}$, $\binom{S_{\mathsf P}}{k-1}$, and $\binom{S_{\mathsf W}}{k - t_k}$ do, where $S_{\mathsf N} > S_{\mathsf P} > S_{\mathsf W}$ are the numbers of primes each must have available.
An implementation with a fixed table of primes invites this picture.
But it ignores the constraint $n \leq B$ (Remark~\ref{rem:binomial}), which the search enforces at every node through the bound $(B/P)^{1/m}$ in line~12, and without that constraint the three counts would differ exponentially.
With it, we do not think they do.
Theorem~\ref{thm:main} gives the same upper bound for all three algorithms, and we expect them to share the leading constant as well (Heuristic~\ref{heur:savings}).
The proof of Theorem~\ref{thm:sharp} suggests why.
The count behaves like a product of binomial coefficients $\binom{N_i}{m_i}$ over ranges of primes near $p_k$.
The smallest primes are nearly forced, and almost all of the freedom lies in the $\Theta(k)$ primes near $p_k$.
The three algorithms differ only in how many of the last choices they leave to a test instead of enumerating: none for $\mathsf N$, one for $\mathsf P$, and $t_k$ for $\mathsf W$.

Algorithm $\mathsf P$ removes only the largest prime from explicit enumeration.
Algorithm $\mathsf W$ removes the largest $t_k$ primes. By
Proposition~\ref{prop:tk}, $t_k<k/2$.
Over our range $t_k$ grows, but $t_k/k$ falls slowly, from about $0.15$ for $k \leq 36$ to $0.11$ at $k = 63$ (Table~\ref{tab:tk}).
This is consistent with $t_k = o(k)$, but the range is too short to establish it.

Let $T_{\mathsf N}(k)$, $T_{\mathsf P}(k)$, and $T_{\mathsf W}(k)$ be the numbers of candidates that the three algorithms examine in a search with $B$ close to $S_k$.

\begin{heuristic}\label{heur:savings}
If $t_k=o(k)$, then $\log T_{\mathsf N}(k)$, $\log T_{\mathsf P}(k)$, and $\log T_{\mathsf W}(k)$ agree up to $o(k)$, and all are at most $(c(\beta_k)+o(1))k$, where $\beta_k$ is defined by $\log S_k=\theta(p_k)+\beta_k k$. The advantage
of $\mathsf W$ over $\mathsf N$ comes from the choices of the largest $t_k$
primes, which contribute $e^{o(k)}$ to the count.
\end{heuristic}

Heuristic~\ref{heur:savings} contains a second assumption, which we make explicit.
Most of the cost of $\mathsf P$ lies in interior nodes, and Theorem~\ref{thm:main} bounds those only by $e^{(c(\beta+1/e)+o(1))k}$.
The heuristic assumes that, at exponential scale, interior nodes are no more numerous than leaves.
This is plausible, since an interior node of depth $\alpha k$ fixes only small primes, which contribute little to the count, and leaves the primes near $p_k$ unchosen, but we have not proved it.

In short, we expect the three algorithms to be equivalent at exponential scale, and the improvement from the
$\mathsf W$ base case to be subexponential, governed by $t_k$. If $t_k$ could be
shown to grow linearly in $k$ (which $t_k<k/2$ does not exclude), $\mathsf W$
would improve the leading constant itself.
The size of $\lambda(P)$ at the leaves controls $t_k$, through $R < \lambda(P)$.
If we had better heuristics for the size of $\lambda(P)$ when $P$ is a product of many small admissible primes, we could estimate $t_k$, and with it the advantage of $\mathsf W$, directly.
We leave this as an open problem.
See \cite{Pollack26} for a step in this direction, an upper bound for the number of $n \leq x$ with $\lambda(n) \leq y$.

\subsection{Admissibility}\label{sec:admissibility}

The tree $\Tree(B,k)$ contains only admissible products.
The counts in Theorems~\ref{thm:main} and~\ref{thm:sharp} do not use this, and are over all square-free products.
How much smaller is the tree?
We do not know in general.
Remark~\ref{rem:adm} shows that the upper bounds still hold and that, for $\beta > 1 + \log 2$, the tree is still exponentially large.
Proposition~\ref{prop:tax} shows that a fixed small prime prunes the tree heavily.
A product containing it must avoid many other primes, and so it must be larger by a factor exponential in $k$.

\begin{remark}[Admissibility]\label{rem:adm}
We record what admissibility does and does not change.
\begin{enumerate}
\item[(a)] The upper bounds are unaffected.
Admissible products are a subset of all square-free products, so Theorem~\ref{thm:main} and Corollary~\ref{cor:runtime} hold as stated.
\item[(b)] Admissibility does not reduce the order of the tree.
If $p < q$ are primes in an interval $(Q/2, Q]$, then $q \not\equiv 1 \pmod p$, since $q \equiv 1 \pmod p$ would force $q \geq 2p + 1 > Q$.
So any set of primes in such an interval is pairwise admissible, and every product of $k$ of them is a depth-$k$ node of $\Tree(B,k)$ once it is at most $B$.
Fix $\varepsilon > 0$ and take $Q = 2(1+\varepsilon)k\log k$.
By the prime number theorem the interval $(Q/2, Q]$ contains $(1+\varepsilon+o(1))k$ primes, and any $k$ of them have product at most $Q^k$, where
\[
  \log Q^k = \theta(p_k) + \bigl(1 + \log 2 + \log(1+\varepsilon) + o(1)\bigr)k .
\]
Hence, for every $\beta > 1 + \log 2$, choosing $\varepsilon$ with $\log(1+\varepsilon) < \beta - 1 - \log 2$ shows that $\Tree(B,k)$ has at least $\binom{(1+\varepsilon+o(1))k}{k} = e^{\Omega(k)}$ nodes of depth $k$.
Together with (a), the tree has size $e^{\Theta(k)}$ in this range.
\item[(c)] We do not know how admissibility affects the constant $c(\beta)$.
The density of primes admissible to $P$ depends on which small primes divide $P$, so it is not a uniform thinning of the primes, and Theorem~\ref{thm:sharp} is a statement about all square-free products.
\end{enumerate}
\end{remark}

Admissibility matters most for the small primes.
A prime $p$ in $P$ makes every prime $q \equiv 1 \pmod p$ inadmissible, which is a proportion $1/(p-1)$ of all primes.
For $p = 3$ this is half of them.
For the search, this is pruning, and the subtree below $p$ is thin.
For the products themselves, it means that a product containing $p$ must be larger, as the next result shows.

\begin{proposition}[Small primes force larger products]\label{prop:tax}
Let $D$ be a fixed finite set of odd primes, and put
\[
  \delta(D) = \prod_{p \in D} \frac{p-2}{p-1} .
\]
As $k \to \infty$, every depth-$k$ node $n$ of $\Tree(B,k)$ that is divisible by every prime in $D$ satisfies
\[
  \log n \geq \theta(p_k) + k \log\bigl(1/\delta(D)\bigr) + o(k) .
\]
In particular, if a fixed prime $p$ divides $S_k$, then $\beta_k \geq \log\frac{p-1}{p-2} - o(1)$.
\end{proposition}

\begin{proof}
Let $M = \prod_{p \in D} p$, and let $\mathcal{Q}_D$ be the set of primes $q \nmid M$ with $q \not\equiv 1 \pmod p$ for every $p \in D$.
Let $q$ be a prime factor of $n$ not in $D$, and let $p \in D$.
If $q > p$, then $q \not\equiv 1 \pmod p$ because $n$ is admissible; if $q < p$, then $q \equiv 1 \pmod p$ is impossible.
So every prime factor of $n$ outside $D$ lies in $\mathcal{Q}_D$.
The set $\mathcal{Q}_D$ is the set of primes in $\prod_{p \in D}(p-2)$ of the $\varphi(M) = \prod_{p \in D}(p-1)$ reduced residue classes modulo $M$.
By the prime number theorem for arithmetic progressions, it contains $(\delta(D) + o(1))\pi(x)$ primes up to $x$.
Hence its $m$-th element $q_m$ satisfies $\log q_m = \log m + \log\log m + \log(1/\delta(D)) + o(1)$.
Since $n$ has at least $k - |D|$ prime factors in $\mathcal{Q}_D$,
\[
  \log n \geq \sum_{m \leq k - |D|} \log q_m = k\log k + k\log\log k - k + k\log\bigl(1/\delta(D)\bigr) + o(k) ,
\]
and the first three terms are $\theta(p_k) + o(k)$.
For the last claim, $S_k$ is a depth-$k$ node of $\Tree(S_k, k)$.
\end{proof}

So a fixed small prime $p$ adds about $k\log\frac{p-1}{p-2}$ to $\log n$.
This is $k\log 2$ for $p = 3$, $k\log(4/3)$ for $p = 5$, and roughly $k/(p-2)$ in general.
Including $p$ saves only the ratio of $p$ to the prime it displaces, which is $e^{O(\log k)}$.
The proposition alone does not exclude any prime from $S_k$ in our range, since $\beta_k \geq 0.80$ there and $\log 2 \approx 0.693$.
But it does show that if $\beta^* < \log\frac{p-1}{p-2}$, then $p \nmid S_k$ for all large $k$.

\subsection{Theory and practice}\label{sec:practice}

The bounds above are far from what we observe.
For the new values, $\beta_k$ falls from about $0.88$ at $k = 40$ to about $0.81$ at $k = 64$ (Table~\ref{tab:sk}).
At $\beta = 0.81$ we have $c(\beta) \approx 2.49$, so Theorem~\ref{thm:main} allows on the order of $e^{2.49k}$ nodes of depth $k$.
At a given $k$ we can also evaluate the bound of Lemma~\ref{lem:rankin} directly.
For $B = 1.05\,S_k$, minimizing it numerically over $\sigma$ and $z$ gives about $10^{41}$ at $k = 40$ and $10^{51}$ at $k = 50$.
In our clean passes with the same $B$, the number of leaves with $P\lambda(P) > B$ was about $3.0 \cdot 10^{8}$ at $k = 40$ and $1.0 \cdot 10^{11}$ at $k = 50$.
That is growth like $e^{0.58k}$, and the counts are smaller than the bound by more than thirty orders of magnitude.

The data do support Conjecture~\ref{conj:U}.
For every $40 \leq k \leq 64$, $S_k$ is smaller than $U(k)$, by a factor of about $10^{18}$ at $k = 40$, growing to about $10^{30}$ at $k = 64$.

By Theorem~\ref{thm:sharp}, Rankin's method is not what loses here.
The bound counts all square-free products, while the tree contains only admissible ones (Remark~\ref{rem:adm}(c)).
It also counts nodes of depth $k$, while the leaves of $\mathsf W$ lie about $t_k$ levels higher.
We have not measured these two effects separately.
Proposition~\ref{prop:tax} shows that each small prime in $P$ prunes the tree heavily, and $P$ contains many small primes.
The lower bound of Remark~\ref{rem:adm}(b) does not reach our data either, since it needs $\beta > 1 + \log 2 \approx 1.69$.
So the exponential order $e^{\Theta(k)}$ is proved only for larger $\beta$ than we observe, and the constant that governs our computations is not explained by the analysis.
Modeling admissibility well enough to predict the observed constant is, we think, the main open problem left by this paper.
We return to this in Section~\ref{sec:lpf}.
There we measure the effect of admissibility directly, by comparing $S_k$ with the least cyclic number with $k$ prime factors, and we look at one visible effect of it, the least prime factor of $S_k$.

\section{Implementation}\label{sec:implementation}

Algorithm~\ref{alg:W} is short, and our implementation of it is not.
In this section we describe the implementation by referring to the lines of Algorithm~\ref{alg:W}.
None of the changes below alter the search tree; they change only how quickly we walk it.
Our code is written in C++ and uses GMP for multiprecision arithmetic.
It is available at \url{https://github.com/jewebste/small-carmichael-numbers}.

\subsection{The preproduct and $\lambda(P)$ (lines 2 and 13)}

Most nodes of the tree are interior, so the cost of line~13 is paid far more often than the cost of a leaf.
We avoid multiprecision arithmetic there whenever possible.

We store $P$ as a GMP integer times a $64$-bit accumulator.
Appending $q$ multiplies only the accumulator, and the accumulator is folded into the GMP integer when the product would overflow.
The same is done for $\lambda(P)$.
We also keep $\lambda(P)$ in factored form, as the exponents of $2$, $3$, $5$, and $7$ together with a short sorted list of larger prime powers.
The point is to keep Euclid's algorithm out of line~13.
Computed directly, $\lcm(\lambda(P), q-1) = \lambda(P)(q-1)/\gcd(\lambda(P), q-1)$ needs a gcd at every interior node, and for the larger $k$ the number $\lambda(P)$ does not fit in a machine word.
Instead we precompute the factorization of $q-1$ for every prime $q$ in range.
Line~13 is then a merge of two short sorted lists, keeping the larger exponent of each prime, and $\lambda(P)$ is multiplied only by the prime powers that grow.
The merge is faster than Euclid's algorithm.
Together with the additive update of $\log\lambda(P)$ described below, it saved $6$ to $11$ percent of the total running time for $20 \leq k \leq 34$.
Often $(q-1) \mid \lambda(P)$, and then the child shares the parent's $\lambda$ and nothing is copied.

Euclid's algorithm still appears once per leaf, in line~3, where $r^\star = P^{-1} \bmod \lambda(P)$ is computed by a modular inversion.
We do not know how to avoid it.
Can $r^\star$ be computed incrementally from the parent, or amortized over the leaves below a common ancestor?
The difficulty is that the modulus changes from $\lambda(P)$ to $\lambda(Pq)$ at each step, so the parent's inverse is an inverse modulo the wrong number.
We leave this as an open question.

We do not store the prime factors of $P$.
A node carries $P$, $\lambda(P)$ in factored form, and the bit vector of admissible primes described below, but no list of the primes in $P$.
The search never needs one.
Admissibility is tracked by the bit vector, and $\lambda(P)$ is updated from the factorization of $q - 1$ alone.
We need the factorization of $n = PR$ only when line~6 accepts it, and that happens a handful of times in a whole search.
At that point the query has already split $R$ into its primes, and we recover the primes of $P$ by dividing $P$ by the primes up to $\pmx(P)$.

We perform the test $P\lambda(P) > B$ in line~2 with floating-point logarithms.
Both $\log P$ and $\log \lambda(P)$ are updated additively in line~13, using a table of $\log q$.
A rounding error here can only move a node to the wrong branch of line~2, which affects the cost but not the correctness (see Section~\ref{sec:algorithm}).
The constant $2$ in line~2 is a parameter of our code, and $2$ is the value we used in practice.
For the larger $k$ in Table~\ref{tab:tk}, constants as large as $5$ could be used with no major change to the behavior of the algorithm.
Replacing $2$ by $c$ affects only nodes that reach $c$ or fewer remaining primes while $P\lambda(P) \leq B$.
But almost every leaf with $P\lambda(P) > B$ has at least $5$ primes remaining; those with $4$ or fewer make up less than $1\%$ of the leaves for each $k \geq 56$ that we measured.
For small $k$, where $t_k$ is closer to $3$, this is no longer true.

\subsection{Admissible primes as bit masks (line 12)}

The loop in line~12 is the innermost loop of the search, and we implement it with bit masks.
The integers up to a compile-time bound $\texttt{MAX\_PRIME}$ that are coprime to $30$ are encoded as bits, eight per block of $30$.
At start-up we clear the composite bits with a sieve of Eratosthenes, and for each prime $q$ we precompute a mask $M_q$ that clears every $r \equiv 1 \pmod q$.
Each node carries its own bit vector of primes admissible to $P$.
In line~13, the child's vector is the parent's vector ANDed with $M_q$, and only the words above $q$ need to be touched.
Line~12 then reduces to repeatedly finding the next set bit (a count-trailing-zeros instruction) until the prime exceeds $(B/P)^{1/m}$.
We compute that bound once per node, with one call to $\exp$, and recompute it only when $B$ changes.

Two cautions are in order.
First, the encoding omits $3$ and $5$, so the search does not see Carmichael numbers divisible by $3$ or $5$.
Nothing is lost in our range.
Every Carmichael number is cyclic, and a branch-and-bound search like the one for $C_k$ in Section~\ref{sec:lpf} shows that, for every $26 \leq k \leq 64$, no cyclic number with $k$ prime factors that is divisible by $3$ or by $5$ is less than $S_k$.
For $3$ the margin is wide: the least cyclic number with $k$ prime factors divisible by $3$ exceeds $S_k$ by a factor of more than $10^{7}$ for $k \geq 40$.
Second, if $(B/P)^{1/m}$ ever exceeds $\texttt{MAX\_PRIME}$, then the search aborts rather than skip primes.
We then recompile with a larger bound and resume at the failed subtree.

\subsection{The pair step (lines 12 and 13)}

Near the root, line~12 runs over small primes, and nearly every child is again an interior node.
For $q \leq 500$ we skip the intermediate level and append two primes at once.
At start-up we tabulate, for each such $q_1$, the primes $q_2 > q_1$ admissible to $q_1$, together with the combined mask $M_{q_1} \wedge M_{q_2}$ and the factorization of $\lcm(q_1 - 1, q_2 - 1)$.
A node then creates its grandchildren $Pq_1q_2$ directly from this table.
This saves the cost of building the intermediate frames.
It does not change the leaves or the work done there.

\subsection{The leaf (lines 3--10)}

In line~4 we start the progression higher than $\pmx(P)$.
The cofactor $R$ is a product of $m$ primes, each larger than $p$, so we may start at the least term exceeding $p^m$.
When the search is run in passes with increasing bounds (see Section~\ref{sec:passes}), we also start above $B'/P$, where $B'$ is the bound of the previous pass.

If $P\lambda(P) > B$, then the loop in lines 5--10 runs at most once.
Indeed, any completion satisfies $0 < R \leq B/P < \lambda(P)$ and $R \equiv r^\star \pmod{\lambda(P)}$, so $R = r^\star$.
The only candidate is $Pr^\star$, and it is a candidate only if $r^\star > p$ and $Pr^\star \leq B$; the next term already has $P(r^\star + \lambda(P)) > P\lambda(P) > B$.
We check $\log P + \log r^\star \leq \log B$ in floating point before forming $Pr^\star$, and we discard $r^\star$ unless $\gcd(r^\star, P) = 1$.

In every computation reported here, this is the only kind of leaf that mattered.
The other branch of line~2, a leaf with $m \leq 2$ and $P\lambda(P) \leq B$, is there for correctness.
We never saw it taken in any meaningful computational regime, and we exercised it only in regimes constructed for testing.
For completeness, our code handles it by sieving the progression: for each prime $\ell \leq \texttt{MAX\_PRIME}$ that is at most $p$ or is not admissible to $P$, we mark the terms divisible by $\ell$, and only the unmarked terms reach line~6.

\subsection{The Carmichael query (line 6)}

We answer line~6 with the Fermat-test and GCD method of \cite{SW25}, restricted to $R$.
Write $n - 1 = 2^s d$ with $d$ odd.
Since $(p - 1) \mid (n - 1)$ for every $p \mid P$ by construction, $b^n \equiv b \pmod P$ holds automatically.
So we need only check $b^n \equiv b \pmod R$, and all of the arithmetic can be done modulo $R$ instead of modulo $n$.
If $m > 1$ and $R$ is prime, we reject at once.
Otherwise, for $b = 2, 3, 5, \ldots$ we compute $b^{d}, b^{2d}, \ldots, b^{2^{s-1}d} \bmod R$, reject if $b^n \not\equiv b \pmod R$, and split each unfactored divisor $D$ of $R$ by $\gcd(D, b^{2^i d} + 1)$.
We stop when every part is prime, and accept if there are exactly $m$ parts $q$, each with $(q - 1) \mid (n - 1)$.

One might ask ``Why not trial divide $R$ first?''
We tried this, and it does not help.
The exponentiations that build confidence that $n$ is a Carmichael number (and so warrant spending computational effort on factoring $R$) produce the factors of $R$ as a byproduct.
By the time enough bases have passed, $R$ has usually been factored already.
A separate pre-filter either rediscovers those factors or is spent on a candidate that would have failed anyway.

\subsection{From recursion to an explicit stack}

The recursion in line~13 is not deep (at most $k$ levels), so the hardware stack is not a concern.
But a recursive call cannot be handed to another thread.
We therefore replaced the hardware stack by a software one.
A frame holds $P$, $\lambda(P)$, the bit vector of admissible primes, the bound $(B/P)^{1/m}$, and the position of the loop in line~12.
The main loop pops a frame, advances its loop by one prime, and pushes the child.
A frame is then a self-contained description of an unexplored subtree, and that is exactly what a parallel scheduler needs.

\subsection{Work stealing}

Our first parallel version cut the tree at a fixed depth and dealt the subtrees out round-robin.
This failed in the familiar way: a few subtrees rooted at small primes contain most of the work, and most threads sat idle while those few finished.

We now use work stealing \cite{BL99}.
Each thread owns a local stack and a deque.
When a thread creates a child frame in line~13, it decides whether the child might be expensive.
If the child's interval in line~12 is wide, in the sense that
\[
  \log\bigl((B/P)^{1/m}\bigr) - \log p > \tau
\]
for a parameter $\tau$, the child is pushed onto the thread's deque; otherwise it stays on the local stack.
A thread takes work from the back of its own deque, and an idle thread steals from the front of another's, so that thieves take the oldest and typically largest subtrees \cite{BL99,CL05}.
We protect each deque with a mutex rather than using a lock-free deque, since the deques see little traffic.
As far as we know, this interval-width test is our own heuristic, and we do not have a citation for it.
We typically used $\tau = 1/2$, chosen after a fair amount of experimentation.
We did not study the choice of $\tau$ systematically, and we leave it as an under-investigated question.

The bound $B$ is shared.
Threads read $\log B$ without locking (a stale value is safe, just slightly wasteful), and a thread that finds a Carmichael number lowers $B$ under a lock.

A master thread enumerates the first four primes of $P$ in order and launches the workers on each resulting subtree.
When the workers rejoin, that subtree is finished, which gives us a checkpoint; a restart skips every subtree up to the last checkpoint.

\subsection{Passes and clean passes}\label{sec:passes}

We do not know $S_k$ in advance, so we run the doubling passes of Definition~\ref{def:doubling}, except that the initial bound $B_0$ comes from a heuristic formula (see the discussion after Corollary~\ref{cor:runtime}).
We take $B_0 = 2^m$, where $m$ is the nearest integer to $\bigl(k\log k + (0.0025k + 1.024)k\bigr)/\log 2$.
The constants are an empirical fit and have no theoretical basis; we compare this formula with Pinch's in Section~\ref{sec:results}.
The program also accepts a starting bound on the command line, in place of the formula.
With the lower bound in line~4, the pass with bound $2B'$ examines only the $n \in (B', 2B']$ that the previous pass did not.
Within a pass, when a Carmichael number is found, we may lower $B$ to it.
To certify $S_k$ we instead run a \textit{clean pass}, with $B$ fixed for the whole pass.
A clean pass is then an exhaustive search of a known range, and its claim is easy to state: no Carmichael number with $k$ prime factors up to $B$ is smaller than the least one it reports.
(The break in line~7 can skip later terms of a progression, but those are larger than the term just reported.)
It also makes the pass independent of the order in which threads happen to find things.
The tree that is searched depends only on $k$ and $B$, so counts such as $t_k$ in Table~\ref{tab:tk} are well defined and repeatable.

\section{Results}\label{sec:results}

We found $S_k$ for every $40 \leq k \leq 64$.
Table~\ref{tab:sk} summarizes the new values, and the complete values are in OEIS A006931 \cite{OEIS}.
Multiplying out each factorization and checking Korselt's criterion shows that each value is a Carmichael number with $k$ prime factors.
Minimality rests on the computation.
For each value we ran a clean pass (Section~\ref{sec:passes}) with $B \geq S_k$, and it found no smaller Carmichael number with $k$ prime factors.
Our code is public, so that anyone can repeat these runs.

\begin{table}[ht]
\centering
\begin{tabular}{rrrrr@{\qquad}rrrrr}
\toprule
$k$ & digits & $p_{\min}$ & $p_{\max}$ & $\beta_k$ & $k$ & digits & $p_{\min}$ & $p_{\max}$ & $\beta_k$ \\
\midrule
40 & 84 & 17 & 1409  & 0.881 & 53 & 118 & 19 & 1297 & 0.838 \\
41 & 87 & 17 & 769   & 0.891 & 54 & 121 & 19 & 1051 & 0.845 \\
42 & 89 & 17 & 1153  & 0.879 & 55 & 124 & 19 & 953  & 0.839 \\
43 & 92 & 17 & 641   & 0.869 & 56 & 126 & 19 & 1009 & 0.832 \\
44 & 95 & 17 & 859   & 0.892 & 57 & 129 & 19 & 1123 & 0.827 \\
45 & 97 & 17 & 601   & 0.873 & 58 & 132 & 19 & 1123 & 0.817 \\
46 & 100 & 17 & 757  & 0.872 & 59 & 134 & 19 & 1321 & 0.817 \\
47 & 102 & 17 & 691  & 0.860 & 60 & 137 & 19 & 1873 & 0.817 \\
48 & 105 & 17 & 1009 & 0.849 & 61 & 140 & 19 & 2689 & 0.801 \\
49 & 108 & 17 & 1289 & 0.864 & 62 & 143 & 19 & 1471 & 0.806 \\
50 & 110 & 17 & 937  & 0.856 & 63 & 145 & 19 & 1201 & 0.807 \\
51 & 113 & 19 & 691  & 0.845 & 64 & 148 & 19 & 1579 & 0.807 \\
52 & 116 & 19 & 1051 & 0.857 &    &     &    &      &       \\
\bottomrule
\end{tabular}
\caption{The new values of $S_k$: the number of decimal digits, the least and greatest prime factors, and $\beta_k$, defined by $\log S_k = \theta(p_k) + \beta_k k$.}
\label{tab:sk}
\end{table}

Pinch suggested that $\log S_k$ might be well approximated by $2\log 2\,(k-1)\log(k-1)$, and noted some support for this from the ratios for $k \leq 35$ \cite{Pinch18}.
Those ratios fall from $1.27$ at $k = 10$ to $0.984$ at $k = 35$, and our values continue the trend, reaching $0.940$ at $k = 64$.
If $\beta^*$ is finite, then the ratio tends to $1/(2\log 2) \approx 0.72$, so the two heuristics cannot both be right.
The value of $\beta_k$ falls from $1.23$ at $k = 10$ to $0.91$ at $k = 35$, and to $0.81$ at $k = 64$ (Table~\ref{tab:sk}), which we take as modest support for the hypothesis of Corollary~\ref{cor:runtime}.
The formula we used for the starting bound (Section~\ref{sec:passes}) is a fit to the same data and has no better claim to the limit.
It implies that $\beta_k$ grows, like $2.02 + 0.0025k - \log\log k$, but it stays within $3.3$ bits of $\log_2 S_k$ for $10 \leq k \leq 64$.

\subsection{Hardware and running times}

We developed and timed the code on a MacBook, using a single core.
Table~\ref{tab:mac} gives the total time to find $S_k$ there, including every doubling pass from the heuristic starting bound.

\begin{table}[ht]
\centering
\begin{tabular}{rr@{\qquad}rr@{\qquad}rr}
\toprule
$k$ & time (s) & $k$ & time (s) & $k$ & time (s) \\
\midrule
20 & 0.12 & 26 & 1.2  & 32 & 36.1 \\
21 & 0.57 & 27 & 9.0  & 33 & 46.9 \\
22 & 0.36 & 28 & 5.9  & 34 & 48.5 \\
23 & 0.34 & 29 & 11.5 & 35 & 200.3 \\
24 & 0.18 & 30 & 13.8 & 37 & 42.6 \\
25 & 2.0  & 31 & 4.7  & 38 & 760.2 \\
   &      &    &      & 39 & 2\,911.1 \\
   &      &    &      & 44 & 44\,033.6 \\
\bottomrule
\end{tabular}
\caption{Time to find $S_k$ on one core of a MacBook, including all doubling passes.}
\label{tab:mac}
\end{table}

The times are far from monotone in $k$.
Most of the variation comes from where the heuristic starting bound lands relative to $S_k$, which decides how many passes are needed and how far the last pass overshoots.

The expensive computations, and every value $S_k$ with $40 \leq k \leq 64$, were run on nodes with an AMD EPYC 9734 processor ($112$ cores), using the work-stealing code of Section~\ref{sec:implementation} with $112$ threads.
For most of these runs we did not keep careful timing data, and the values of $k$ missing from Table~\ref{tab:amd} are those for which we have none.
Table~\ref{tab:amd} gives the times we have, each for a single clean pass with $B$ a little above $S_k$.

\begin{table}[ht]
\centering
\begin{tabular}{rrr@{\qquad}rrr}
\toprule
$k$ & $B/S_k$ & time (s) & $k$ & $B/S_k$ & time (s) \\
\midrule
40 & 1.05 & 1\,012   & 50 & 1.05 & 5\,305 \\
41 & 1.05 & 1\,226   & 51 & 1.05 & 4\,511 \\
42 & 1.05 & 1\,284   & 52 & 1.40 & 55\,138 \\
43 & 1.05 & 1\,378   & 54 & 1.77 & 80\,779 \\
45 & 1.05 & 1\,833   & 55 &      & 75\,297 \\
46 & 1.05 & 2\,222   & 56 &      & 30\,639 \\
47 & 1.05 & 2\,219   & 59 &      & 48\,064 \\
48 & 1.05 & 2\,258   & 60 &      & 111\,646 \\
49 & 1.05 & 5\,312   & 61 &      & ${\approx}\,400\,000$ \\
   &      &          & 63 &      & 712\,670 \\
\bottomrule
\end{tabular}
\caption{Wall-clock time of one clean pass with bound $B$, on one AMD EPYC 9734 node with $112$ threads.}
\label{tab:amd}
\end{table}

These times, too, are not monotone in $k$, for two reasons.
The ratio $B/S_k$ varies, and the work in a pass grows quickly with $B$.
Also, the wall-clock time of a parallel pass depends on how well the work balances, and a single large subtree can dominate.

\section{Least prime factors and cyclic numbers}\label{sec:lpf}

The least prime factor of $S_k$ changes slowly.
It is at least $11$ for every $k \geq 15$, it is $17$ for $31 \leq k \leq 50$, and it is $19$ for $51 \leq k \leq 64$ (Table~\ref{tab:sk}).
By Proposition~\ref{prop:tax}, a small prime $p$ adds about $k\log\frac{p-1}{p-2}$ to $\log S_k$, while it saves only the ratio of the prime it displaces, near the top of the product, to $p$.
The cost grows linearly in $k$ and the saving only logarithmically, so each fixed small prime is eventually excluded.
A crude balance of the two suggests a least prime factor of order $k/\log\log k$.

To see how much of this is admissibility alone, we compare $S_k$ with $C_k$, the least cyclic number with $k$ prime factors, which ignores the Korselt condition.
We computed $C_k$ for every $k \leq 64$ by branch and bound over the same tree as Algorithm~\ref{alg:W}.
That search is easy, for two reasons: taking admissible primes greedily gives a good bound $B$ before the search starts, and at a node $P$ with $m$ primes left, the product of the next $m$ primes admissible to $P$ is a sharp lower bound for every completion.
Our values agree with the terms of OEIS A264907 \cite{OEIS-cyclic}, which are listed for $k \leq 17$.
Table~\ref{tab:cyclic} compares $C_k$ with $S_k$.

\begin{table}[ht]
\centering
\begin{tabular}{rrrrrr}
\toprule
$k$ & $p_{\min}(C_k)$ & $p_{\min}(S_k)$ & $\beta(C_k)$ & $\beta_k$ & $\log_{10}(S_k/C_k)$ \\
\midrule
10 & 5  & 11 & 0.863 & 1.231 & 1.6 \\
15 & 7  & 11 & 0.845 & 1.079 & 1.5 \\
20 & 11 & 11 & 0.831 & 1.013 & 1.6 \\
25 & 11 & 11 & 0.804 & 0.990 & 2.0 \\
30 & 11 & 11 & 0.788 & 0.948 & 2.1 \\
35 & 11 & 17 & 0.767 & 0.913 & 2.2 \\
40 & 11 & 17 & 0.752 & 0.881 & 2.2 \\
45 & 11 & 17 & 0.749 & 0.873 & 2.4 \\
50 & 11 & 17 & 0.738 & 0.856 & 2.6 \\
53 & 11 & 19 & 0.731 & 0.838 & 2.5 \\
54 & 19 & 19 & 0.727 & 0.845 & 2.8 \\
55 & 19 & 19 & 0.724 & 0.839 & 2.7 \\
60 & 19 & 19 & 0.707 & 0.817 & 2.9 \\
64 & 19 & 19 & 0.696 & 0.807 & 3.1 \\
\bottomrule
\end{tabular}
\caption{The least cyclic number $C_k$ and the least Carmichael number $S_k$ with $k$ prime factors. Here $\beta(C_k)$ is defined by $\log C_k = \theta(p_k) + \beta(C_k)k$.}
\label{tab:cyclic}
\end{table}

Three things stand out.
First, $S_k$ is close to $C_k$: the ratio $S_k/C_k$ is at most $10^{3.1}$ for $k \leq 64$.
The part of $\beta_k$ not explained by admissibility, $\beta_k - \beta(C_k)$, falls from $0.37$ at $k = 10$ to $0.11$ at $k = 64$.
So at $k = 64$, admissibility alone accounts for about $86\%$ of $\beta_k$, and the Korselt condition for the rest.
Second, the least prime factor of $C_k$ jumps from $11$ to $19$ at $k = 54$, skipping $13$ and $17$, and from then on it agrees with the least prime factor of $S_k$ for every $k \leq 64$.
Third, for $12 \leq k \leq 53$ the least prime factor of $S_k$ is at least that of $C_k$, and often larger.
We do not think the Korselt condition acts on the small primes directly.
A leaf has a Korselt completion only by chance, so $S_k$ is the least of the comparatively few admissible products below it that happen to have one.
It therefore resembles a typical admissible product of its size more than the extreme one, $C_k$, and we expect a typical product to contain fewer small primes, since each small prime excludes a residue class from the rest of the product (Proposition~\ref{prop:tax}).

\begin{conjecture}\label{conj:lpf}
The least prime factor of $S_k$ tends to infinity with $k$, and it is $o(p_k)$.
Moreover, it equals the least prime factor of $C_k$ for all sufficiently large $k$.
\end{conjecture}

We note a by-product.
Nielsen asked in \cite{OEIS-cyclic} whether the prime factors of $C_k$ always form a normal sequence of primes, that is, whether $C_k$ is obtained by choosing a first prime and then repeatedly taking the least prime admissible to all those already chosen.
They do not.
For $k = 31$, the cyclic number
\[
  11 \cdot 17 \cdot 19 \cdot 29 \cdot 31 \cdot 37 \cdots 197 \cdot 211 \cdot 241 \approx 1.46 \cdot 10^{59}
\]
that we found is already smaller, by a factor of $1.04$, than the product of the first $31$ terms of every normal sequence.
This conclusion does not depend on the minimality of our search.
The greedy choice from $11$ takes $13$, and admitting $13$ excludes $53$, $79$, $131$, and $157$, which turns out to cost more than $13$ saves.
For every $31 \leq k \leq 64$, $C_k$ is not the product of a normal sequence.

\section*{Use of AI tools}
This paper was prepared with substantial assistance from Claude (Anthropic), a large language model, used through Claude Code.
The algorithm of Section~\ref{sec:algorithm} is due to the authors.
The implementation described in Section~\ref{sec:implementation} was written by Claude under our direction, mainly with Claude Sonnet 4.6 and Claude Sonnet 5.
The analysis of Section~\ref{sec:analysis}, including the proofs (among them the proof of Theorem~\ref{thm:sharp}) and Heuristic~\ref{heur:savings}, was carried out primarily by Claude Opus 5.5 under our direction, and was checked by the authors.
Claude Opus 5.5 also drafted much of the text from our outlines and revisions, and helped locate and verify references.
All editorial and stylistic decisions were made by the first author.
The authors take full responsibility for the contents of this paper.

\section*{Acknowledgements}
We thank Andrew Shallue for helpful discussions.
We thank Aidan Johnson for his work on the 2025 project and for researching numerical estimators for $S_k$.
Thanks also to Frank Levinson for generous donations to Butler University that helped make this computation possible.


\begin{thebibliography}{99}
\bibitem{Korselt1899} A.~Korselt, Probl\`eme chinois, \emph{L'Interm\'ediaire
  des Math\'ematiciens} \textbf{6} (1899), 142--143.
\bibitem{Simerka1885} V.~\v{S}imerka, Zbytky z arithmetick\'e posloupnosti,
  \emph{\v{C}asopis pro p\v{e}stov\'an\'i mathematiky a fysiky} \textbf{14}
  (1885), 221--225.
\bibitem{Carmichael1910} R.~D.~Carmichael, Note on a new number theory
  function, \emph{Bull.\ Amer.\ Math.\ Soc.}\ \textbf{16} (1910), 232--238.
\bibitem{Poulet1938} P.~Poulet, Tables des nombres compos\'es v\'erifiant le
  th\'eor\`eme de Fermat pour le module~2 jusqu'\`a 100.000.000,
  \emph{Sphinx} (Brussels) \textbf{8} (1938), 42--45.
\bibitem{Swift1975} J.~D.~Swift, Table of Carmichael numbers to $10^9$,
  \emph{Math.\ Comp.}\ \textbf{29} (1975), 338--339.
\bibitem{PSW80} C.~Pomerance, J.~L.~Selfridge, and S.~S.~Wagstaff, Jr., The
  pseudoprimes to $25\cdot10^9$, \emph{Math.\ Comp.}\ \textbf{35} (1980),
  1003--1026.
\bibitem{Jaeschke1990} G.~Jaeschke, The Carmichael numbers to $10^{12}$,
  \emph{Math.\ Comp.}\ \textbf{55} (1990), 383--389.
\bibitem{Keller1988} W.~Keller, The Carmichael numbers to $10^{13}$,
  \emph{Abstracts Amer.\ Math.\ Soc.}\ \textbf{9} (1988), 328--329.
\bibitem{Guthmann1992} A.~Guthmann, On the computation of Carmichael numbers,
  Tech.\ Report 218, Universit\"at Kaiserslautern, 1992.
\bibitem{Pinch93} R.~G.~E.~Pinch, The Carmichael numbers up to $10^{15}$,
  \emph{Math.\ Comp.}\ \textbf{61} (1993), 381--391.
\bibitem{Pinch17} R.~G.~E.~Pinch, The Carmichael numbers up to $10^{17}$,
  arXiv:math/0504119 (2005).
\bibitem{Pinch18} R.~G.~E.~Pinch, The Carmichael numbers up to $10^{18}$,
  arXiv:math/0604376 (2006).
\bibitem{AGHS14} W.~R.~Alford, J.~Grantham, S.~Hayman, and A.~Shallue,
  Constructing Carmichael numbers through improved subset-product algorithms,
  \emph{Math.\ Comp.}\ \textbf{83} (2014), 899--915.
\bibitem{PinchTable} R.~G.~E.~Pinch, Tables of Carmichael numbers,
  \url{http://www.s369624816.websitehome.co.uk/rgep/cartable.html}.
\bibitem{Erdos56} P.~Erd\H{o}s, On pseudoprimes and Carmichael numbers,
  \emph{Publ.\ Math.\ Debrecen} \textbf{4} (1956), 201--206.
\bibitem{AGP94} W.~R.~Alford, A.~Granville, and C.~Pomerance, There are
  infinitely many Carmichael numbers, \emph{Ann.\ of Math.}\ \textbf{139}
  (1994), 703--722.
\bibitem{MRC25} A.~Johnson, S.~Webster, and J.~Webster, poster, Butler MRC
  2025, \url{https://blue.butler.edu/~jewebste/MRC_2025.pdf}.
\bibitem{OEIS} The On-Line Encyclopedia of Integer Sequences, Sequence
  A006931, \url{https://oeis.org/A006931}.
\bibitem{OEIS-cyclic} The On-Line Encyclopedia of Integer Sequences, Sequence
  A264907, \url{https://oeis.org/A264907}.
\bibitem{SW25} A.~Shallue and J.~Webster, Algorithms for Carmichael numbers,
  arXiv:2506.09903 (2025).
\bibitem{AKS} M.~Agrawal, N.~Kayal, and N.~Saxena, PRIMES is in P,
  \emph{Ann.\ of Math.}\ \textbf{160} (2004), 781--793.
\bibitem{Pollack26} P.~Pollack, Small values of Carmichael's lambda
  function, arXiv:2606.19132 (2026).
\bibitem{Rankin38} R.~A.~Rankin, The difference between consecutive prime
  numbers, \emph{J.\ London Math.\ Soc.}\ \textbf{13} (1938), 242--247.
\bibitem{EN81} P.~Erd\H{o}s and J.-L.~Nicolas, Sur la fonction ``nombre de
  facteurs premiers de $n$'', \emph{Enseign.\ Math.}\ (2) \textbf{27} (1981),
  3--27.
\bibitem{Pomerance89} C.~Pomerance, Two methods in elementary analytic number
  theory, in R.~A.~Mollin (ed.), \emph{Number Theory and Applications},
  Kluwer, Dordrecht, 1989, 135--161.
\bibitem{HT88} A.~Hildebrand and G.~Tenenbaum, On the number of prime factors
  of an integer, \emph{Duke Math.\ J.}\ \textbf{56} (1988), 471--501.
\bibitem{HT86} A.~Hildebrand and G.~Tenenbaum, On integers free of large
  prime factors, \emph{Trans.\ Amer.\ Math.\ Soc.}\ \textbf{296} (1986),
  265--290.
\bibitem{Tenenbaum} G.~Tenenbaum, \emph{Introduction to Analytic and
  Probabilistic Number Theory}, 3rd ed., Graduate Studies in Mathematics
  \textbf{163}, Amer.\ Math.\ Soc., 2015.
\bibitem{HR18} G.~H.~Hardy and S.~Ramanujan, Asymptotic formulae in
  combinatory analysis, \emph{Proc.\ London Math.\ Soc.}\ (2) \textbf{17}
  (1918), 75--115.
\bibitem{BL99} R.~D.~Blumofe and C.~E.~Leiserson, Scheduling multithreaded
  computations by work stealing, \emph{J.\ ACM} \textbf{46} (1999), 720--748.
\bibitem{CL05} D.~Chase and Y.~Lev, Dynamic circular work-stealing deque,
  \emph{Proc.\ SPAA} 2005, 21--28.
\end{thebibliography}
\end{document}